\documentclass[final,3p,times]{elsarticle}
\usepackage{amsmath,hyperref}
\usepackage{lineno}

\usepackage{epstopdf}
\journal{XXX}

\usepackage{lineno,hyperref}
\usepackage{graphicx}% Include figure files
\usepackage{dcolumn}% Align table columns on decimal point
\usepackage{bm}% bold math
\usepackage{epsfig}
\usepackage{booktabs}
\usepackage{subfigure}
\usepackage{graphics}
\usepackage{amssymb}
\usepackage{amsmath}
\usepackage{array}
\usepackage{color}
\usepackage{booktabs}
\usepackage{multirow}
\usepackage{caption}
\usepackage{chngpage}
\usepackage{subfigure}
\usepackage{mathrsfs}

\biboptions{numbers,sort&compress}
\modulolinenumbers[1]

\begin{document}

\begin{frontmatter}

\title{A hybrid physics-informed neural network for nonlinear partial differential equation}
\author[mymainaddress,mysecondaddress]{Chunyue Lv}
\author[mymainaddress,mysecondaddress]{Lei Wang\corref{mycorrespondingauthor}}
\cortext[mycorrespondingauthor]{Corresponding author}
\ead{leiwang@cug.edu.cn; wangleir1989@126.com}
\author[mymainaddress,mythirdaddress]{Chenming Xie}

\address[mymainaddress]{School of Mathematics and Physics, China University of Geosciences, Wuhan 430074, China}
\address[mysecondaddress]{Center for Mathematical Sciences, China University of Geosciences, Wuhan 430074, China}
\address[mythirdaddress]{Kindo medical data technology company limited , Wuhan 430073, China}

\begin{abstract}
The recently developed physics-informed machine learning  has made great progress for solving nonlinear partial differential equations (PDEs), however, it may fail to provide reasonable approximations to the PDEs with discontinuous solutions. In this paper, we focus on the discrete time physics-informed neural network (PINN), and propose a hybrid PINN scheme for the nonlinear PDEs. In this approach, the local solution structures are classified as smooth and non-smooth scales by introducing a discontinuity indicator, and then the automatic differentiation technique is employed for resolving smooth scales, while an improved weighted essentially non-oscillatory (WENO) scheme is adopted to capture discontinuities. We then test the present approach by considering the viscous and inviscid Burgers equations , and it is shown that compared with original discrete time PINN, the present hybrid approach has a better performance in approximating the discontinuous solution even at a relatively larger time step.

\end{abstract}


\end{frontmatter}

%\linenumbers

\section{Introduction}
Nonlinear partial differential equations (PDEs) have been the focus of many studies due to their frequent appearance in many applications in physical, biological, and social sciences \cite{{birk1997non}}. Because of most nonlinear PDEs do not have analytical solutions expect for some extremely simple situations, and therefore, various numerical methods have been developed during the past several decades, such as the finite-difference \cite{mazumder2015}, finite-volume
\cite{mazumder2015} and finite-element \cite{lazarov1996}, and the lattice Boltzmann methods \cite{shi2009pre}, to name but a few.

Apart from the traditional numerical approaches, the rise of deep learning in recent years offers an opportunity to develop surrogate model for solving nonlinear PDEs \cite{dissana1994,sirigan2018,wu2020data,raissi2019physics}. Among various neural network based PDEs solvers, the physics-informed neural network (PINN) proposed by Raissi et al. \cite{raissi2019physics} has drawn considerable attention due to its effectiveness in solving both forward and inverse nonlinear PDE problems as well as the straightforward implementation \cite{raissi2019physics,karniadakis2021}.
The core idea of PINN algorithm is to encode the underlying physical laws into the loss function of the neural network such that the governing equations are enforced by minimizing the residual loss function with  automatic differentiation \cite{raissi2019physics,karniadakis2021}. More recently, the PINN and its variants have been successfully applied to solve various PDEs, including the Navier-Stokes equation, the Kdv equations \cite{raissi2019physics}, heat transfer equations \cite{zobeiry2021a}, Fokker-Plank equation \cite{xu2021solving} and so on. Despite the great progress of PINNs in solving nonlinear PDEs, the governing equations appeared in previous works are usually the so-called diffusion-dominated equations, and it is known that the solutions in these cases are usually smooth enough such that the traditional PINNs are capable of achieving good predictive accuracy according to the classical universal approximation theorem \cite{cybenko1989appro}. However, as pointed out by some authors, it is a challenge for current PINNs to predict the solutions of PDEs with advection-dominant character, in which the solutions often develop discontinuities in some finite time even if the initial conditions is smooth enough \cite{dafermo2005hyper}. To address this issue, some enhanced PINNs are proposed recently.  By using more scattered points around the discontinuous regions, Mao et al. \cite{mao2020physics} investigated the possibility of continuous time PINN in solving the high-speed aerodynamic flows, and found that compared with random or uniform training points, the predicted solutions with clustered training points are more efficient in some cases.  In addition, another work reported by Fuks et al. \cite{fuks2020limitions} shows that the performance of PINN in solving hyperbolic PDEs can be enhanced by adding a small diffusion term to the governing PDEs. More recently, by combining the traditional finite volume methodology with PINNs, Patel et al. \cite{patel2021therm} proposed a thermodynamically consistent PINNs for hyperbolic systems, and the solutions of various PDEs like Euler, Buckley-Leverett equations are well predicted with this method. Although there are some works on utilization of PINNs for solving the PDEs with discontinuities, the neural network appeared in some previous works either needs a priori knowledge of the solution  \cite{mao2020physics}  or the performance  of the PINNs largely depends on the neural network architecture \cite{fuks2020limitions,almajid2022pre}.

In this work, we propose a hybrid PINN (hPINN) for nonlinear PDEs, and the basic idea of this approach is to introduce a discontinuity indicator into the neural network such that the local solution structures can be classified as smooth and non-smooth scales. Then the automatic differentiation technique (one of the most important techniques in deep learning field \cite{baydin2018auto}) is employed for resolving smooth scales, while the improved weighted essentially non-oscillatory (WENO) scheme \cite{jiang1996eff}, which is a famous shock-capturing scheme, is adopted to capture discontinuities. The reminder of the present paper is organized as follows. In Section 2, the details of the proposed approach are presented. Then, some results for viscous and inviscid Burgers equations are presented in Section 3, and a brief conclusion is drawn in Section 4.

\section{Methodology}
In this work, we are interested in designing a hybrid PINN model for the following nonlinear PDE,
\begin{equation}
{u_t} + f{\left( u \right)_{\emph{\textbf{x}}}} = \upsilon g{\left( u \right)_{\emph{\textbf{xx}}}} + h\left( {\emph{\textbf{x}},t} \right),
\label{eq1}
\end{equation}
in which $u$ is the scalar variable and usually a function of both time $t$ and space $\emph{\textbf{x}}$, $h\left( {\emph{\textbf{x}},t} \right)$ is the source term, and $\upsilon$ is the diffusion coefficient. In addition, the first derivations in Eq. (\ref{eq1}) are associated with convection while the second derivatives are responsible for diffusion.

\subsection{Hybrid physics-informed neural network}
\begin{figure}[h]
\centering
\includegraphics[width=0.9\textwidth]{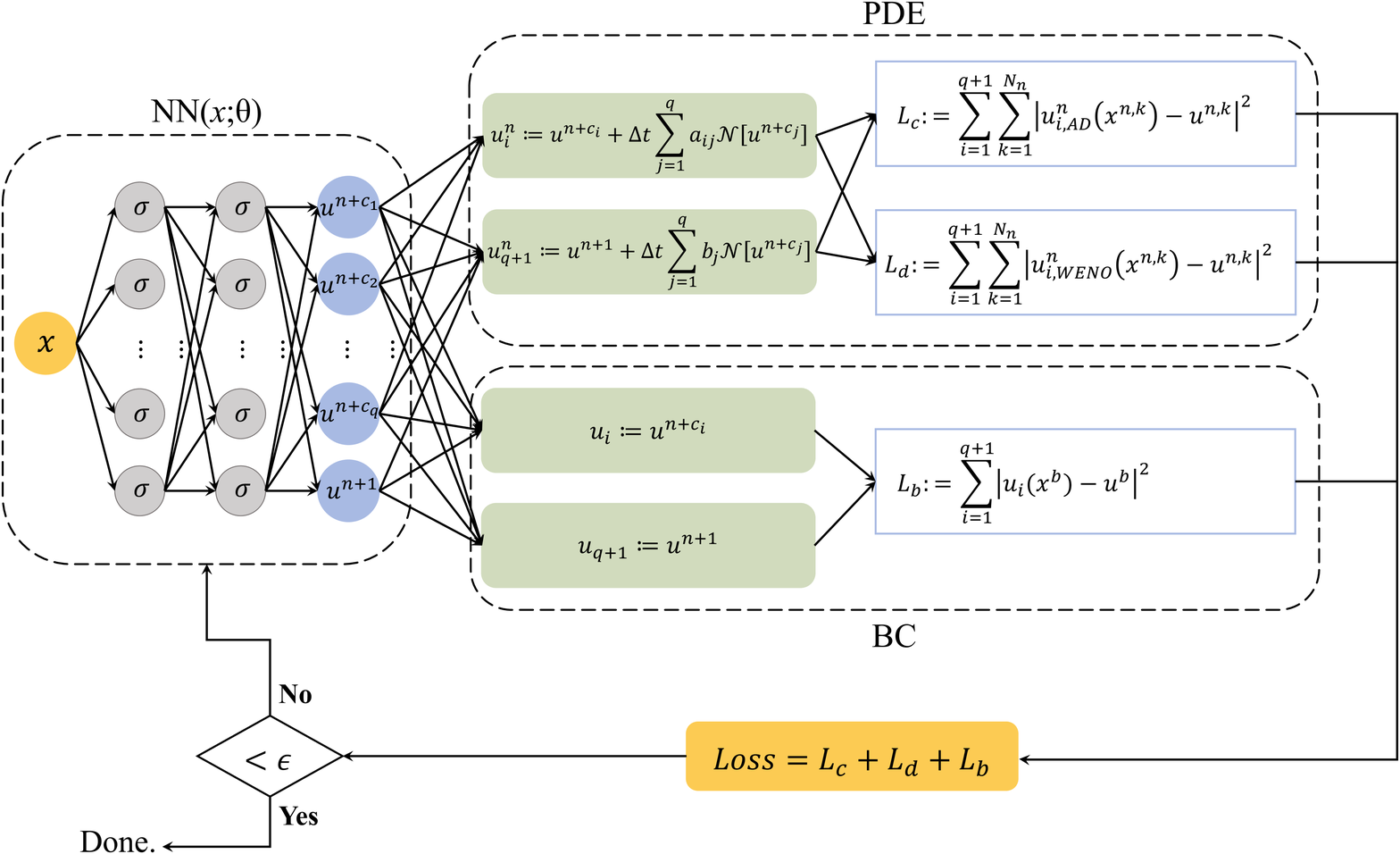}
\caption{Schematic of the hPINN for solving the nonlinear PDEs, in which the spatial vector ${x}$ is the input, $\sigma$ represents the activation function, which is set to be the hyperbolic tangent activation function in this work \cite{raissi2019physics}, and all the hyperparameters in the neural network are initialized by the Xavier initialization method \cite{meng2020pp}. } \label{fig1}
\end{figure}
As presented by Raissi et al. \cite{raissi2019physics}, the generally used PINNs can be classified as the continuous time model and the discrete time model. Compared with the latter model, the continuous time model has drawn tremendous interest in the PINN community \cite{mao2020physics,fuks2020limitions,patel2021therm,almajid2022pre}. However, as pointed out by Raissi et al. \cite{raissi2019physics}, in order to enforce the physical conservation laws , the continuous time model needs a large amount of collocation points in the entire spatio-temporal domain, and thus the training of the original continuous time model  is usually prohibitively expensive. Fortunately, this limitation could be well addressed by employing the discrete time model \cite{raissi2019physics}. In such a case, the present work intends to solve the forward PDE problems by using the discrete time model, and the model proposed here can be viewed as the  variant of the original discrete time model.

To have a better understand on the present approach, the original discrete time PINN is first reviewed. The main idea of this model is to discrete the the time operator in Eq. (\ref{eq1}) by using the classical Runge-Kutta methods with $q$ stage \cite{raissi2019physics}, then we can obtain the following equations
\begin{equation}
\begin{split}
& {u^{n + {c_i}}}\left( x \right) = {u^n}\left( x \right) - \Delta t\sum\limits_{j = 1}^q {{a_{ij}}{\mathscr N}} \left[ {{u^{n + {c_j}}}\left( x \right)} \right],\;i = 1, \ldots ,q, \\
& {u^{n + 1}}\left( x \right) = {u^n}\left( x \right) - \Delta t\sum\limits_{j = 1}^q {{b_j}{\mathscr N}} \left[ {{u^{n + {c_j}}}\left( x \right)} \right],
\end{split}
\end{equation}
in which  ${\mathscr N}\left[ u \right] = f{\left( u \right)_x} - \upsilon g{\left( u \right)_{xx}} - h\left( x \right)$, ${u^{n + {c_i}}}\left( x \right) = u\left( {{t^n} + {c_i}\Delta t,x} \right)$, and the above time discretization scheme may be implicit or explicit, depending on the values of $\left\{ a_{ij},b_j,c_j \right\}$. In such a case, the PINN model usually consists of two sub-networks, in which the first neural network ${{{NN}}(w,b)}$ takes  $\emph{\textbf{x}}$ as input and outputs the values of $\emph{\textbf{u}}$ at different time stages, i.e., $\left[ {{u^{n + {c_1}}}, \ldots ,{u^{n + {c_q}}},{u^{n + 1}}} \right]$, and they are further fed into the second network to encode the governing equations as well as the corresponding boundary conditions. For the forward problems considered here, since the initial/boundary conditions as well as the governing equations value are always known, the loss function in the original discrete time PINN model is given by \cite{raissi2019physics},

\begin{equation}
\begin{split}
&L = {L_{PDE}} + {L_{BC}},\\
&{L_{PDE}} = {\sum\limits_{j = 1}^q {\sum\limits_{i = 1}^{{N_n}} {\left| {u_j^n\left( {{x^{n,i}}} \right) - {u^{n,i}}} \right|} }^2} + {\left| {u_{q + 1}^n\left( {{x^{n,i}}} \right) - {u^{n,i}}} \right|^2},\\
&{L_{BC}} = \sum\limits_{j = 1}^q {{{\left| {u_j^n\left( {{x^{b}}} \right) - {u^{b}}} \right|}^2}}  + {\left| {u_{q + 1}^n\left( {{x^{b}}} \right) - {u^{b}}} \right|^2},
\end{split}
\end{equation}
where $N_n$ is the number of the collocation points, which can be randomly or uniformly sampled inside the computational domain,  $\left\{ {{x^{n,i}},{u^{n,i}}} \right\}_{i = 1}^{{N_n}}$ is the numerical data obtained at time-step $t^n$, $x^b$ represents the collocation points on the boundary, and $u_i^n$ is defined as
\begin{equation}
\begin{split}
& u_i^n = {u^{n + {c_i}}} + \Delta t\sum\limits_{j = 1}^q {{a_{ij}}{\mathscr N}} \left[ {{u^{n + {c_j}}}\left( x \right)} \right], \;i = 1, \ldots ,q, \\
& u_{q + 1}^n = {u^{n + 1}} + \Delta t\sum\limits_{j = 1}^q {{b_j}{\mathscr N}} \left[ {{u^{n + {c_j}}}\left( x \right)} \right].
\end{split}
\end{equation}

As presented by some previous works \cite{fuks2020limitions,almajid2022pre, patel2021therm}, although the original PINN works well for various nonlinear PDEs, it may fail to provide reasonable approximations to the PDEs with dominant convection.  To address this issue, an intuitive way is to employ some  discontinuity-capturing methods in the neural network. In this work, inspired by the weighted essentially nonoscillatory schemes (WENO) \cite{jiang1996eff,jiang1995jcp}, which are a popular class of numerical methods in capturing discontinuities, an improved WENO scheme \cite{jiang1996eff}, dubbed WENO-Z, is incorporate into the PINN. On the basis of this idea, a simple way is to discretize the differential operator in the whole computational domain by using WENO-Z approach. However, this treatment may lead to another two issues: For one thing, when applying the WENO-Z scheme in the whole domain, the training cost of the neural network will be increased significantly due to the local characteristic decomposition and the nonlinear-weights computing \cite{jiang1996eff,jiang1995jcp}. For another, the advantage of the automatic differentiation technique used in original PINN will be weakened. In this setting, an alternative solution is to introduce the hybrid concept, for which the efficient automatic differentiation technique is adopted to compute all the differential operators of the PDEs in smooth regions while the expensive WENO-Z scheme is just applied in the vicinity of discontinuities, and it is observed that this idea is widely employed in scientific computing community \cite{fu2019a}.

Fig. \ref{fig1} shows a schematic of the present \emph{hybrid} PINN (hPINN). Note that the evolution of the discontinuities depends little on the viscous term (i.e., $g(u)_{xx}$), the automatic differentiation technique is adopted to compute the second derivative in both the smooth and non-smooth regions for its efficient.  For the approximation of the first order derivative ( i.e., $f(u)_x$) in the non-smooth scale, the fifth-order WENO-Z scheme is used. For brevity in the presentation, we use one-dimensional scalar nonlinear PDEs as an example, and assume the scatter points {$x_j$} are uniform. The cell size as well as the cells are denoted by $\Delta x = {x_{j + 1}} - {x_j}$ and ${I_j} = \left[ {{x_{j - 1/2}},{x_{j + 1/2}}} \right]$ with ${x_{j + 1/2}} = {x_j} + \Delta x/2$, respectively. Then use a conservative finite difference scheme to the first-order derivative term in Eq. (\ref{eq1}), thus obtain
\begin{equation}
\frac{{\partial f}}{{\partial x}}\left| {_{x = {x_i}}} \right. = \frac{{{{\hat f}_{i + 1/2}} - {{\hat f}_{i - 1/2}}}}{{\Delta x}},
\end{equation}
in which the numerical fluxes ${{\hat f}_{i + 1/2}}$ and ${{\hat f}_{i - 1/2}}$ are the positive and negative parts of $f(u)$ at the cell boundaries.  To ensure the numerical stability, the global Lax-Friedrichs flux splitting is used for calculating the flux $f(u)$
\begin{equation}
f\left( u \right) = {f^ + }\left( u \right) + {f^ - }\left( u \right),
\end{equation}
where ${f^ \pm }\left( u \right) = \frac{1}{2}\left( {f\left( u \right) \pm \lambda u} \right)$ and $\lambda  = \max \left| {f'\left( u \right)} \right|$, and with this, we have ${{\hat f}_{i + 1/2}} = {{\hat f}^ + }_{i + 1/2} + {{\hat f}^ - }_{i + 1/2}$. Since the negative part of the flux ${{\hat f}_{i - 1/2}}$ is symmetric to the positive part with respect to $x_{i+\frac{1}{2}}$, here we will only describe how ${{\hat f}^ + }_{i + 1/2}$ is constructed, and for convenience, we will drop the "+" sign in the superscript. As presented in \cite{jiang1996eff}, the fifth-order WENO flux  ${{\hat f}_{i + 1/2}}$ is defined as
\begin{equation}
{{\hat f}_{i + 1/2}} = {\omega _0}{{\hat f}^{\left( 0 \right)}}_{i + 1/2} + {\omega _1}{{\hat f}^{\left( 1 \right)}}_{i + 1/2} + {\omega _2}{{\hat f}^{\left( 2 \right)}}_{i + 1/2},
\end{equation}
in which ${{\hat f}^{\left( k \right)}}_{i + 1/2}\;\left( {k = 0,1,2} \right)$ are the flux values assigned to three stencils $\left\{ {{x_{j - 2}},{x_{j - 1}},{x_j}} \right\}$, $\left\{ {{x_{j - 1}},{x_{j }},{x_{j+1}}} \right\}$ and $\left\{ {{x_{j}},{x_{j + 1}},{x_{j+2}}} \right\}$,  and they are given by
\begin{equation}
\begin{split}
& {{\hat f}^{\left( 0 \right)}}_{i + 1/2} = \frac{1}{6}\left( {2{f_{j - 2}} - 7{f_{j - 1}} + 11{f_j}} \right) ,\\
& {{\hat f}^{\left( 1 \right)}}_{i + 1/2} = \frac{1}{6}\left( { - {f_{j - 1}} + 5{f_j} + 2{f_{j + 1}}} \right), \\
& {{\hat f}^{\left( 2 \right)}}_{i + 1/2} = \frac{1}{6}\left( {2{f_j} + 5{f_{j + 1}} - {f_{j + 2}}} \right).\\
\end{split}
\end{equation}
In addition, the nonlinear weights ${\omega _0}$, ${\omega _1}$  and ${\omega _2}$ are calculated as follow \cite{jiang1996eff}
\begin{equation}
{\omega _k} = \frac{{{\alpha _k}}}{{{\alpha _0} + {\alpha _1} + {\alpha _2}}},\;{\alpha _k} = {d_k}\left[ {1 + {{\left( {\frac{{{\tau _5}}}{{{\beta _k} + \varepsilon }}} \right)}^2}} \right],k = 0,1,2,
\end{equation}
where $\varepsilon $ is a small positive sensitivity parameter, which is used to avoid divisions by zero in the weights formulation, and is chosen to be $10^{-40}$ in this work as suggested in \cite{jiang1996eff}. Coefficient $d_{k}$ is the optimal weights given by
\begin{equation}
{d_0} = \frac{1}{{10}},{d_1} = \frac{6}{{10}},{d_2} = \frac{3}{{10}}.
\end{equation}
$\tau$ a global smoothness indicator used in WENO-Z scheme \cite{jiang1996eff}, which is characterized by
\begin{equation}
{\tau _5} = \left| {{\beta _0} - {\beta _2}} \right|,
\end{equation}
In addition, ${{\beta _k}}$ is the smoothness indicator, which measures the smoothness of the solution over a particular stencil, and are defined as
\begin{equation}
\begin{split}
&{\beta _0} = \frac{{13}}{{12}}\left( {{f_{j - 2}} - 2{f_{j - 1}} + {f_j}} \right) + \frac{1}{4}{\left( {{f_{j - 2}} - 4{f_{j - 1}} + 3{f_j}} \right)^2}, \\
&{\beta _1} = \frac{{13}}{{12}}\left( {{f_{j - 1}} - 2{f_j} + {f_{j + 1}}} \right) + \frac{1}{4}{\left( {{f_{j + 1}} - {f_{j - 1}}} \right)^2}, \\
&{\beta _2} = \frac{{13}}{{12}}\left( {{f_j} - 2{f_{j + 1}} + {f_{j + 2}}} \right) + \frac{1}{4}{\left( {3{f_j} - 4{f_{j + 1}} + {f_{j + 2}}} \right)^2}. \\
\end{split}
\end{equation}

Finally, in order to distinguish the non-smooth scales from the smooth regions in the neural network, an adaptive discontinuity indicator is needed. Following the work of Fu \cite{fu2019a}, the separation between discontinuities and smooth regions can be achieved by
\begin{equation}
{\gamma _k} = \frac{1}{{{{\left( {{\beta _k} + \delta } \right)}^p}}},\;k = 0,1,2,3,
\end{equation}
where $\delta$ is equal to $10^{-4}$ and $10^{-3}$ for one-dimensional and multi-dimensional problems,  $q$ is a scale separation parameter, the value of it is set to 6, as suggested in \cite{fu2019a}, ${\beta _3}$ is another smoothness indicator over the stencil of $\left\{ {{x_{j + 1,}}{x_{j + 2,}}{x_{j + 3}}} \right\}$, and it is given by
\begin{equation}
{\beta _3} = \frac{1}{3}\left( {{f_{i + 1}}\left( {22{f_{i + 1}} - 73{f_{i + 2}} + 29{f_{i + 3}}} \right) + {f_{i + 2}}\left( {61{f_{i + 2}} - 49{f_{i + 3}}} \right) + 10{f_{i + 3}}{f_{i + 3}}} \right).
\end{equation}
Then the discontinuity indicator is defined as
\begin{equation}
{\rm{Flag = }}\left\{ \begin{array}{l}
0,\;{\chi _k} > {C_T},\;\forall k \in \left\{ {0,1,2,3} \right\}\\
1,\;{\rm{otherwise}}
\end{array} \right.,
\end{equation}
in which $\chi _k$ the normalized smoothness indicator given by ${\chi _k} = {{{\gamma _k}} \mathord{\left/{\vphantom {{{\gamma _k}} {\sum\limits_{k = 0}^3 {{\gamma _k}} }}} \right. \kern-\nulldelimiterspace} {\sum\limits_{k = 0}^3 {{\gamma _k}} }}\;\left( {k = 0, \ldots ,3} \right)\;$ and $C_{T}$ is a built-in parameter defined as ${C_T} = 5 \times {10^{ - 4}}$.

\section{Computational results and discussion}
\begin{figure}[h]
\centering
 \subfigure[]{ \label{fig2a}
\includegraphics[width=0.48\textwidth]{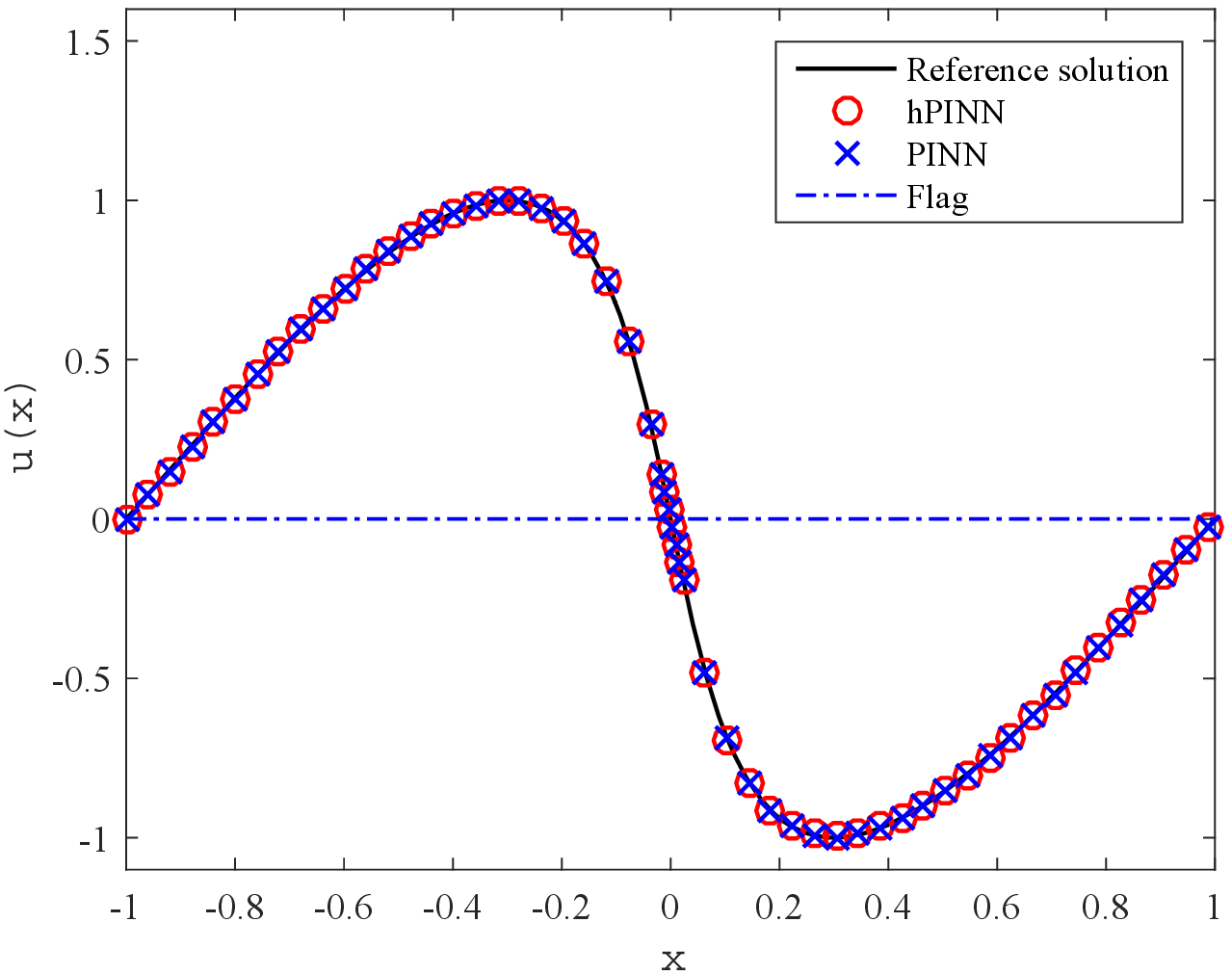}}
 \subfigure[]{ \label{fig2b}
\includegraphics[width=0.48\textwidth]{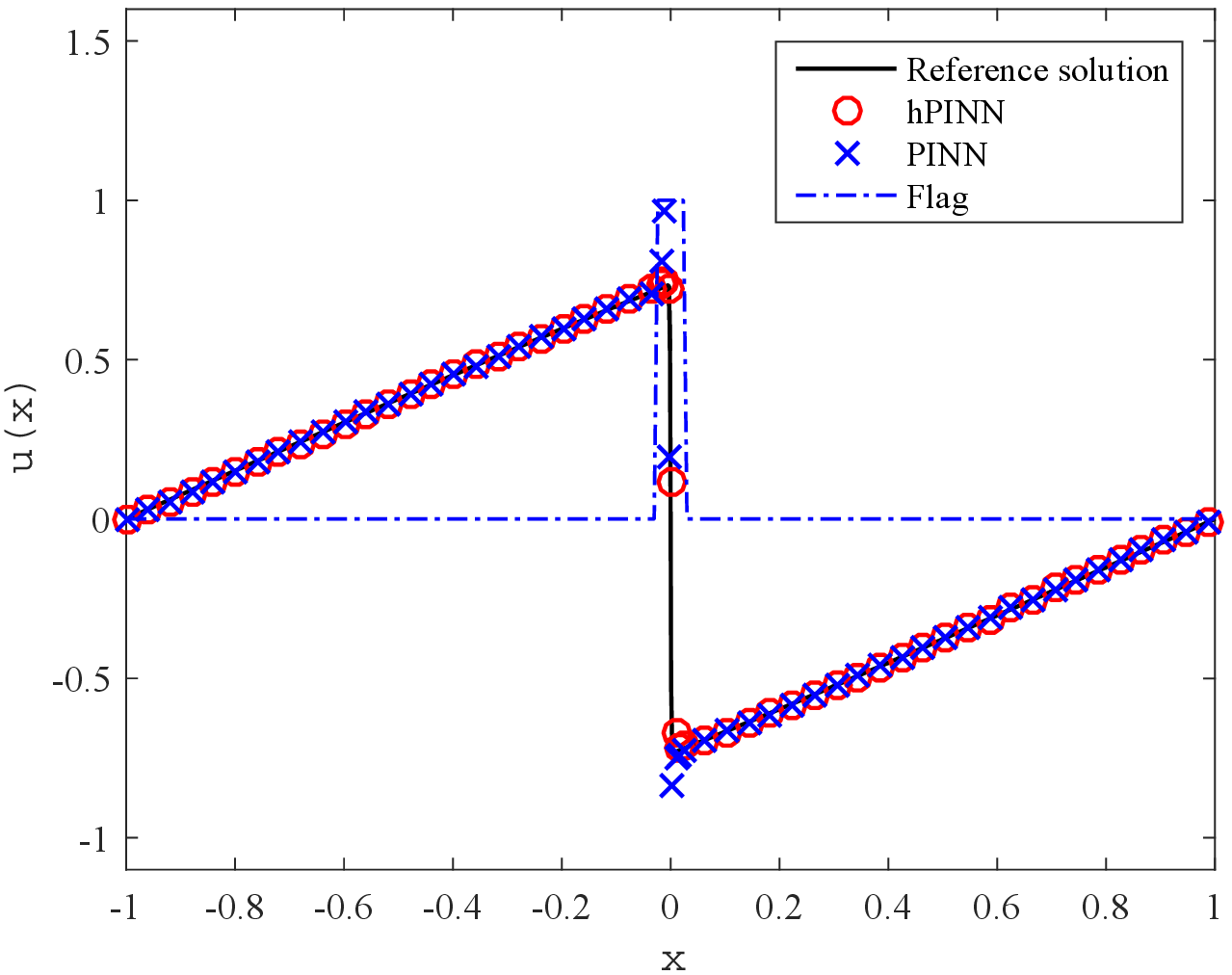}}
\caption{Comparison of the viscous Burgers equation at $t=0.2$ (a) and $t=1.0$ (b), in which $\nu=10^{-4}/\pi$, and the reference solution is computed with WENO-Z scheme using 1000 mesh points.}
 \label{fig2}
\end{figure}
In this section, the one-dimensional viscous and inviscid Burgers equations with Dirichlet boundary conditions are chosen to validate the present hPINN,
\begin{equation}
\begin{split}
&\frac{{\partial u}}{{\partial t}} + \frac{{\partial \left( {\frac{{{u^2}}}{2}} \right)}}{{\partial x}} = \nu \frac{{{\partial ^2}u}}{{\partial {x^2}}},\;x \in \left[ { - 1,1} \right],\;t \in \left[ {0,1} \right],\\
&u\left( {0,x} \right) =  - \sin \left( {\pi x} \right),\\
&u\left( {t, - 1} \right) = u\left( {t,1} \right) = 0,
\end{split}
\end{equation}
and the predicted results are compared with the original PINN. As for the neural networks architecture, we employ 5 hidden-layers with 20 neuron in each layer and the number of residual points are 300 in training process, and the learning rate is set as $10^{-4}$. In both cases, the network training is performed up to the loss function is smaller than $10^{-5}$.

\begin{figure}[h]
\centering
 \subfigure[]{ \label{fig3a}
\includegraphics[width=0.48\textwidth]{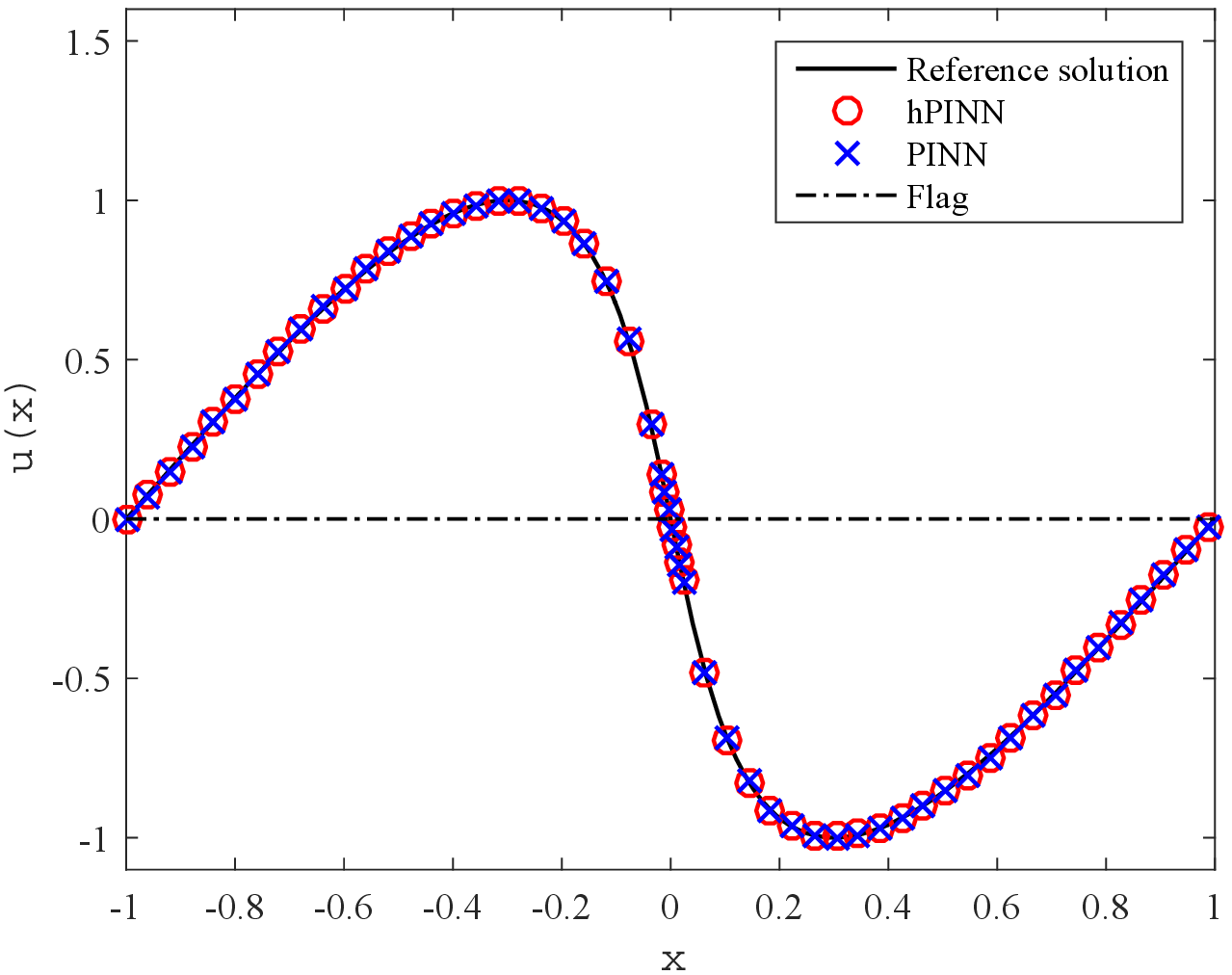}}
 \subfigure[]{ \label{fig3b}
\includegraphics[width=0.48\textwidth]{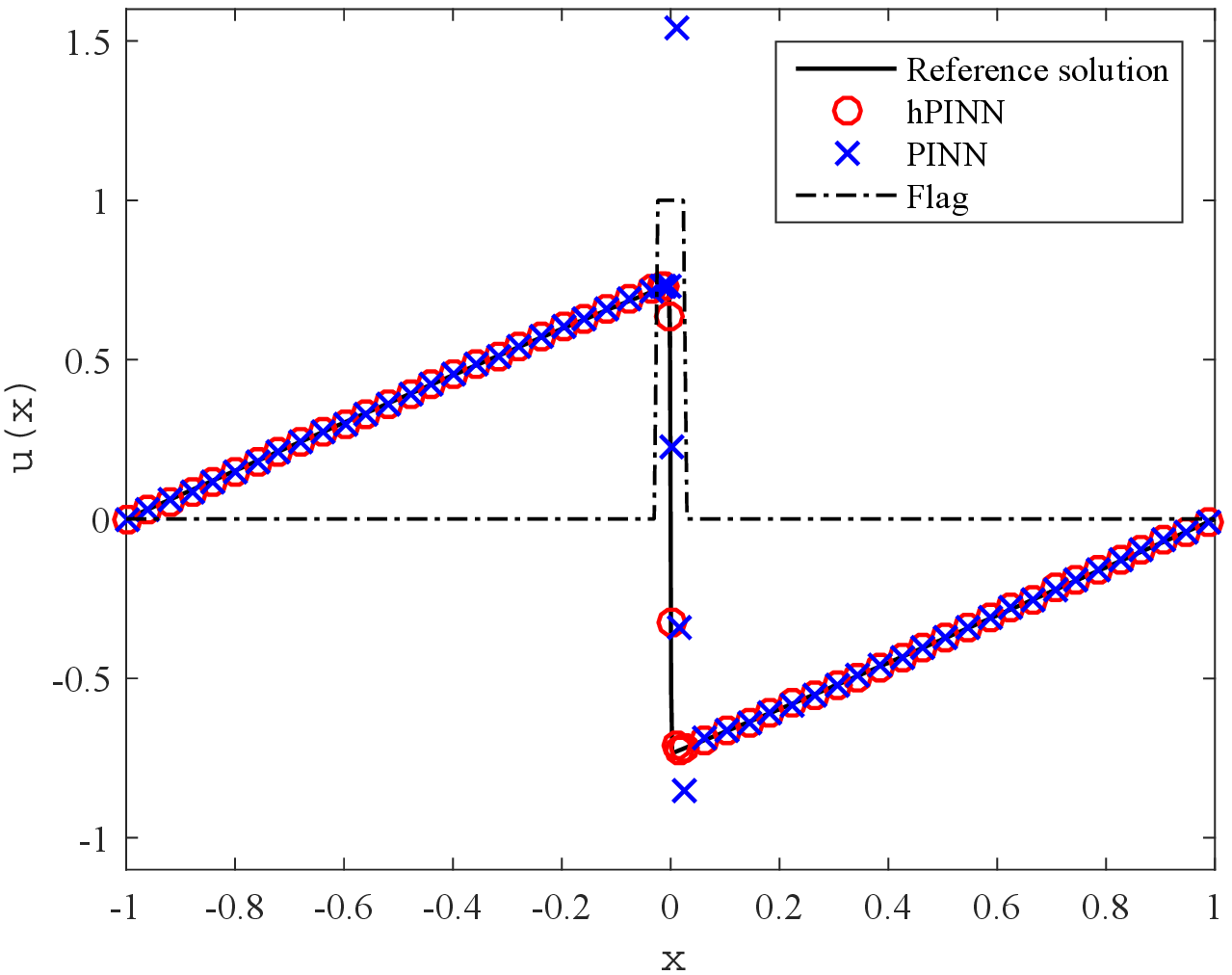}}
\caption{Comparison of the inviscid Burgers equation (i.e., $\nu=0.0$) at $t=0.2$ (a) and $t=1.0$ (b), in which the reference solution is computed with WENO-Z scheme using 1000 mesh points.}
 \label{fig3}
\end{figure}

Fig. \ref{fig2} and Fig. \ref{fig3} illustrate the distributions of the predict solutions for both the viscous and inviscid Burgers equations, in which the reference solutions are obtained by using fifth order WENO-Z scheme with 1000 mesh points. As depicted in these figures, it can be clearly seen that in regions where the solution is smooth both the PINN and hPINN work well, when the shock wave appears in the solution, the original PINN fails to provide reasonable approximations at the shock position, and this phenomenon is more distinct for the inviscid Burgers equation (see Fig. \ref{fig2b} and Fig. \ref{fig3b}).  Different from the original PINN, the results obtained with present hPINN converges to the reference solution with a non-oscillatory profile even in the inviscid case, which further suggests that the performance of present hPINN is better that previous PINN.

To demonstrate the robustness of the present hPINN, the global relative errors obtained at various Runge-Kutta stages $q$ and the time-steps $\Delta t$ for present hPINN are also investigated, and the corresponding results are shown in Table \ref{Tab3-1}, in which the results obtained by using original PINN are not incorporate since it fails to predict reasonable results. As shown in this table, one can find that for a fixed Runge-Kutta stage $q$/time-step $\Delta t$, the global relative error is always increased/decresed with increasing  time-step $\Delta t$/ Runge-Kutta stage $q$. Moreover,  due to the distinct feature of the neural network, the Runge-Kutta stage used here can be an arbitrarily larger number, making it possible to  solve the nonlinear PDEs with a very large time step, which is largely different from the traditional numerical method.

\begin{table}[h]
\caption{Global relative errors obtained with hPINN for different Runge-Kutta stages $q$ and time-step $\Delta t$ at $t=0.6$.}
 \label{Tab3-1}
 \centering
\begin{tabular}{cccccccc}
\hline
           & \multicolumn{3}{c}{${\nu} = 10^{-4}/\pi$}  &  & \multicolumn{3}{c}{$ \nu = 0.0$}     \\
\cline{2-4} \cline{6-8}
q          & $\Delta t=0.1$       & $\Delta t=0.3$       & $\Delta t=0.6$      & & $\Delta t=0.1$       & $\Delta t=0.3$       & $\Delta t=0.6$          \\
\hline
1          & $2.3438\times10^{-3}$ & $3.4727\times10^{-1}$ & $4.3328\times10^{-1}$ & & $4.4754\times10^{-3}$ & $3.8972\times10^{-2}$ & $5.9224\times10^{-1}$  \\
4          & $2.8555\times10^{-3}$ & $9.2065\times10^{-3}$ & $3.3564\times10^{-1}$ & & $2.0912\times10^{-3}$ & $5.3934\times10^{-3}$ & $1.2390\times10^{-1}$  \\
10         & $1.8478\times10^{-3}$ & $2.6818\times10^{-2}$ & $4.2223\times10^{-2}$ & & $3.0013\times10^{-3}$ & $5.1584\times10^{-3}$ & $9.8035\times10^{-3}$  \\
50         & $2.2558\times10^{-3}$ & $3.5475\times10^{-3}$ & $9.6065\times10^{-2}$ & & $2.3808\times10^{-3}$ & $3.7173\times10^{-3}$ & $1.1018\times10^{-2}$  \\
 \hline
\end{tabular}
\end{table}

\section{Conclusion}
\label{sec4}
In this paper, by incorporating a discontinuity indicator into the neural network, a hybrid physics-informed neural network (hPINN) is proposed, in which the automatic differentiation technique is employ to compute all the differential operators of the PDEs in smooth regions, while the classical WENO-Z scheme is adopted to compute the convection term in the vicinity of discontinuities. Then the viscous and inviscid Burgers are selected to validated the present model. Based on the results, we observe that the present hPINN is able to provide reasonable approximation even at the shock position, and overall performance of the present hPINN is better than the original PINN. Further,  due to the distinct feature of the neural network, the Runge-Kutta stage used here can be an arbitrarily larger number, making it possible to  solve the nonlinear PDEs with a very large time step, which is very hard or even impossible in traditional numerical method. Finally, we would like to point out that although the present work focuses on the one-dimensional problems, it is easily to extent the present model to the multi-dimensional problems, which will be considered in ongoing work.

\section*{CRediT authorship contribution statement}

\textbf{Chunyue Lv:} Software, Methodology, Investigation. \textbf{Lei Wang:} Conceptualization, Methodology, Writing -
review  ${\rm{\& }}$  editing.

\section*{Acknowledgements}
This work is financially supported by the National Natural Science Foundation of China (Grant No. 12002320) and the Fundamental Research Funds for the Central Universities (Grant No. CUGGC05),.

\end{document}